\input amstex
\input epsf
\magnification=\magstep1 
\baselineskip=13pt
\documentstyle{amsppt}
\vsize=8.7truein \CenteredTagsOnSplits \NoRunningHeads
\topmatter
 
\title Computing the probability of intersection \endtitle 
\author Alexander Barvinok  \endauthor
\address Department of Mathematics, University of Michigan, Ann Arbor,
MI 48109-1043, USA \endaddress
\email barvinok$\@$umich.edu \endemail
\date September 13,  2025 \enddate
\thanks  This research was partially supported by NSF Grant DMS 1855428. 
\endthanks 
\abstract  Let $\Omega_1, \ldots, \Omega_m$ be probability spaces, let $\Omega=\Omega_1 \times \cdots \times \Omega_m$ be their product and let 
$A_1, \ldots, A_n \subset \Omega$ be events.
Suppose that each event $A_i$ depends on $r_i$ coordinates of a point  $x \in \Omega$, $x=\left(\xi_1, \ldots, \xi_m\right)$,
and that for each event $A_i$ there are $\Delta_i$ of other events $A_j$ that depend on some of the coordinates that $A_i$ depends on.
Let $\Delta=\max\{5,\ \Delta_i: i=1, \ldots, n\}$ and let $\mu_i=\min\{r_i,\ \Delta_i+1\}$ for $i=1, \ldots, n$.  We prove that if
  ${\bold P}(A_i) < (3\Delta)^{-3\mu_i}$ for all $i$,
   then for any $0 < \epsilon < 1$, the probability ${\bold P}\left( \bigcap_{i=1}^n \overline{A}_i\right)$ of the intersection of the complements of all $A_i$  can be computed within relative error  $\epsilon$ in polynomial time from the probabilities ${\bold P}\left(A_{i_1} \cap \ldots \cap A_{i_k}\right)$ of $k$-wise intersections of the events $A_i$ for 
 $k = e^{O(\Delta)} \ln (n/\epsilon)$.
 \endabstract
\endtopmatter
\document

\head 1. Introduction and main results \endhead 

\subhead (1.1) The setup \endsubhead
Let $\Omega$ be a probability space and let $A_1, \ldots, A_n \subset {\bold \Omega}$ be events. We want to compute (or approximate) 
$${\bold P}\left( \bigcap_{i=1}^n \overline{A}_i \right)=1- {\bold P}\left( \bigcup_{i=1}^n A_i\right), \tag1.1.1$$
 the probability of the intersection of their complements. This is of course a fairly general problem, with many applications. In the usual interpretation, $A_i$ are some ``bad" events and 
 so we are interested in the probability that none of them occurs. 
  If the events $A_1, \ldots, A_n$ are independent, then 
  $${\bold P}\left( \bigcap_{i=1}^n \overline{A}_i\right) = \prod_{i=1}^n \left(1 - {\bold P}(A_i)\right), \tag1.1.2$$
  so to compute the left hand side of (1.1.2) exactly, we only need to know the probabilities ${\bold P}(A_i)$ themselves. 
  
We show that if the dependencies are controlled and the probabilities ${\bold P}(A_i)$ are sufficiently small, then the value of (1.1.1) can be efficiently approximated within an arbitrarily small relative error $0 < \epsilon < 1$ from the probabilities 
  $${\bold P}\left( A_{i_1} \cap \ldots \cap A_{i_k}\right) \quad \text{where} \quad k = O\bigl( \ln (n/\epsilon)\bigr)$$
of the intersections of logarithmically many events $A_i$. 

Let $\Omega_1, \ldots, \Omega_m$ be probability spaces and let ${\boldsymbol \Omega} = \Omega_1 \times \cdots \times \Omega_m$ be their product. 
We write a point $x \in {\boldsymbol \Omega}$ as $x=\left(\xi_1, \ldots, \xi_m\right)$, where $\xi_j \in \Omega_j$, and call $\xi_j$ {\it coordinates} of $x$. We consider complex-valued random variables 
$f: {\boldsymbol \Omega} \longrightarrow {\Bbb C}$. For such an $f$, we define the {\it set of coordinates that $f$ depends on}, as a subset $J_f \subset \{1, \ldots, m\}$ such that for any two $x', x'' \in {\boldsymbol \Omega}$, where $x'=\left(\xi_1', \ldots, \xi_m'\right)$ and $x''=\left(\xi_1'', \ldots, \xi_m''\right)$, we have $f(x')=f(x'')$ whenever 
$\xi_j' = \xi''_j$ for all $j \in J_f$ and there is no subset $J' \subsetneq J_f$ with that property. We say that $f$ {\it depends on at most $r$ coordinates} if $|J_f| \leq r$. We say that two random variables $f, g: {\boldsymbol \Omega} \longrightarrow {\Bbb C}$ {\it share a coordinate} if $J_f \cap J_g \ne \emptyset$, and any coordinate $\xi_j$ with 
$j \in J_f \cap J_g$ is called a {\it shared coordinate} of $f$ and $g$.

Let $A \subset {\boldsymbol \Omega}$ be an event. The {\it indicator} of $A$ is the random variable $[A]: {\boldsymbol \Omega} \longrightarrow \{0, 1\}$ defined by 
$$[A](x)=\cases 1 &\text{if\ } x \in A \\ 0 & \text{if\ } x \notin A. \endcases \tag1.1.3$$
Similarly, we say that $A$ {\it depends on at most $r$ coordinates} if $[A]$ depends on at most $r$ coordinates, and that two events $A, B \subset {\boldsymbol \Omega}$ {\it share a coordinate} if $[A]$ and $[B]$ share a coordinate.

By $\overline{A}$ we denote the complement of an event $A \subset {\boldsymbol \Omega}$, so that $[\overline{A}]=1-[A]$. We also note that $[A] [B] = [A \cap B]$ for events $A, B \subset {\boldsymbol \Omega}$.

We can now state the main result.

\proclaim{(1.2) Theorem} 
Let $\Omega_1, \ldots, \Omega_m$ be probability spaces and let ${\boldsymbol \Omega} = \Omega_1 \times \cdots \times \Omega_m$ be their product. Let 
$A_1, \ldots, A_n \subset {\boldsymbol \Omega}$ be events such that each event $A_i$ depends on at most $r_i$ coordinates and for 
each event $A_i$ at most $\Delta_i$ other events $A_j$ share a coordinate with $A_i$. Let
$$\Delta=\max\left\{5, \ \Delta_i: \ i=1, \ldots, n \right\}$$ and let 
$$\mu_i=\min\{ r_i,\ \Delta_i+1\} \quad \text{for} \quad i=1, \ldots, n.$$
Then for any $0 < \epsilon < 1$ there is $$K=e^{O(\Delta)} \ln (n/\epsilon)$$ such that whenever 
$${\bold P}(A_i)  \ < \ (3 \Delta)^{-3\mu_i}  \quad \text{for} \quad i=1, \ldots, n, \tag1.2.1$$
the value of ${\bold P}\left( \bigcap_{i=1}^n \overline{A}_i \right)$, up to relative error $\epsilon$, is determined by the probabilities 
of $k$-wise intersections
$${\bold P}\left(A_{i_1} \cap \ldots \cap A_{i_k}\right), \quad \text{where} \quad k \leq K. \tag1.2.2$$
More precisely, for $1 \leq k \leq K$, let
$$\sigma_k=\sum_{1 \leq  i_1 < \ldots < i_k \leq n} {\bold P}\left(A_{i_1} \cap \ldots \cap A_{i_k}\right).\tag1.2.4$$
Then for given $\Delta$, $n$ and $0< \epsilon < 1$, there is a polynomial $Q=Q_{\Delta, n, \epsilon}$ in $K$ variables and with $\deg Q \leq K$, such that 
$$\left| \ln {\bold P}\left( \bigcap_{i=1}^n \overline{A}_i\right) - Q\left(\sigma_1, \ldots, \sigma_k\right) \right| \ \leq \ \epsilon,$$
provided (1.2.1) holds.
 
 The implied constant in the ``$O$" notation is absolute.
\endproclaim
Note that an additive $\epsilon$-approximation of $\ln {\bold P}\left( \bigcap_{i=1}^n \overline{A}_i\right)$ translates into a relative $\epsilon$-approximation of 
${\bold P}\left( \bigcap_{i=1}^n \overline{A}_i\right)$. We are interested in the regime where $r_i$ and $\Delta$ are small (fixed in advance), while $n, m$ and $1/\epsilon$ are allowed to grow.

The polynomial $Q$ has rational coefficients, and given $\sigma_1, \ldots, \sigma_K$,  its value $Q(\sigma_1, \ldots, \sigma_K)$ can be computed efficiently, in polynomial in $K$ time, so the main difficulty is in computing $\sigma_k$.
Computing probabilities (1.2.2) reduces to operating with $R=r_{i_1} + \ldots + r_{i_k}$ coordinates and  often can be done roughly in exponential in $R$ time, which leads to an approximation algorithm of a quasi-polynomial 
$n^{O(\ln(n/\epsilon))}$ complexity, if $\Delta$ and $r_i$ are bounded above by constants, fixed in advance. Moreover, in some situations the combinatorial methods of \cite{PR17} and \cite{LSS19} may sharpen it to a genuinely polynomial time algorithm (more on that is in Section 2.2 below).

We note that one can always assume that $r_i \leq 2^{\Delta_i}$ (we will not use it). Indeed, let us choose an event, without loss of generality, $A_n$. Let $I \subset \{1, \ldots, n-1\}$ be the set of indices $i$ such that $A_i$ shares a coordinate with $A_n$, so that $|I| \leq \Delta_n$. Let $J \subset \{1, \ldots, m\}$ be the set of indices $j$ such that $A_n$ depends on the coordinate $\xi_j$, so that $|J| \leq r_n$. With each $j \in J$ we associate a subset $I_j \subset I$ consisting of all indices $i \in I$ such that $A_i$ depends on $\xi_j$. Now, all coordinates $\xi_j$ with the same associated subset $I_j$ can be ``lumped together" into a single coordinate. For example, if 
$A_n$ does not share a coordinate with any other event $A_i$, then all the coordinates that $A_n$ depends on can be considered as a single coordinate in the probability product space, and if $A_n$ shares a coordinate with exactly one other event, say $A_1$, then all the coordinates shared by $A_n$ and $A_1$ can be considered as a single coordinate in the product space, while all remaining coordinates, if any, that $A_n$ depends on can be considered as another coordinate. Since there are $2^{|I|} \leq 2^{\Delta_n}$ subsets of $I$, in the end we obtain at most $2^{\Delta_n}$ coordinates for $A_n$.

The bound (1.2.1) is unlikely to be optimal, though some upper bound in terms of $\Delta$ is definitely needed, since otherwise even deciding whether ${\bold P}\left(\bigcap_{i=1}^n \overline{A}_i\right) >0$ can be quite hard. It would be interesting to find out whether there has to be an exponential dependence on $\mu_i$. 

We now give a simple example of how Theorem 1.2 can be applied.

\example{(1.3) Example: counting integer points} Let $\ell_1, \ldots, \ell_n: {\Bbb R}^m \longrightarrow {\Bbb R}$ be polynomials. In the cube $C=[0, c]^m \subset {\Bbb R}^m$, where $c$ is a positive integer, we consider the set 
$$S =\Bigl\{ x \in C \cap {\Bbb Z^m}: \quad   \ell_i(x) \ \leq \ 0 \quad \text{for} \quad i=1, \ldots, n \Bigr\} \tag1.3.1$$
of integer points, satisfying a system of polynomial inequalities.

Suppose further that each polynomial $\ell_i$ depends on at most $r$ coordinates of $x \in {\Bbb R}^n$, $x=\left(\xi_1, \ldots, \xi_m\right)$, and that for each polynomial
$\ell_i$ there are not more than $\Delta$ of other polynomials $\ell_j$ that depend on some of the same coordinates. We assume that $r$ and $\Delta$ are fixed in advance, while $n$, $m$ and $c$ can vary. We want to approximate the cardinality $|S|$ of $S$.

For $i=1, \ldots, m$, we interpret the set ${\Bbb Z} \cap [0, c]$ of integer points in an interval as a probability space $\Omega_i$ with the uniform measure, and let $\Omega= \Omega_1 \times \cdots \times \Omega_m$. We further introduce events 
$$A_i=\Bigl\{x \in C \cap {\Bbb Z}^m: \quad \ell_i(x) > 0 \Bigr\} \quad \text{for} \quad i=1, \ldots, n$$
and see that 
$$|S| =(c+1)^m {\bold P}\left( \bigcap_{i=1}^n \overline{A}_i\right).$$
Suppose now that 
$${\bold P}(A_i)= |A_i| (c+1)^{-m} \ < \ (3 \Delta)^{-3r} \quad \text{for} \quad i=1, \ldots, n.$$
Theorem 1.2 then asserts that up to a relative error $0 < \epsilon < 1$, the cardinality of $S$  is determined by the cardinalities
$$A_{i_1} \cap \ldots \cap A_{i_k} = \Bigl\{x \in  C\cap {\Bbb Z}^m: \quad \ell_{i_1}(x) > 0, \ldots, \ell_{i_k}(x) >0 \Bigr\} \tag1.3.2$$
for $k=O\left(\ln (n/\epsilon)\right)$. We note that the inequalities in (1.3.2) involve at most $k r$ coordinates $\xi_j$, so that we can compute the probability of each intersection (1.3.2) in the most straightforward way by enumerating $(c+1)^{rk}$ integer points, which in the end results in an approximation algorithm to count points in (1.3.1)
of  $\left(m(c+1)\right)^{O\left( \ln (n/\epsilon)\right)}$ complexity.

If $\deg \ell_i \leq 1$ for $i=1, \ldots, n$, we can do better: in this case, $A_i$ is the set of integer points in a semi-open polyhedron, and the exact enumeration in (1.3.2) can be done with 
$$\bigl((1 + \ln c)(kr)\bigr)^{O\bigl(kr \ln (kr)\bigr)}$$ complexity, see
\cite{Bar08}, which leads to a quasi-polynomial algorithm for approximating the number of integer points in (1.3.1).

Similar algorithms can be designed for computing volumes, in which case the probability spaces $\Omega_i$ are identified with intervals and hence become infinite, see also \cite{GK94} for the complexity of (deterministic) volume computation in small dimensions.
\endexample

\subhead (1.4) Connections \endsubhead The Lov\'asz Local Lemma \cite{EL75} provides a sufficient criterion for the event $\bigcap_{i=1}^n \overline{A}_i$ to have positive probability.
Namely, for each $i=1, \ldots, n$, let $\Gamma_i \subset \{1, \ldots, n\}$ be the set of indices $j \ne i$ such that $A_j$ shares a coordinate with $A_i$. 
Suppose that there are numbers $0 < s_i < 1$ for $i=1, \ldots, n$, such that 
$${\bold P}(A_i) \ \leq \ s_i \prod_{j \in \Gamma_i} (1-s_j) \quad \text{for} \quad i=1, \ldots, n. \tag1.4.1$$
Then 
$${\bold P}\left(\bigcap_{i=1}^n \overline{A}_i \right) \ \geq \ \prod_{j=1}^n (1-s_i) \ > \ 0. \tag1.4.2$$
As in Theorem 1.2, let $\Delta$ be an upper bound on the number of events $A_j$ that share a coordinate with any particular event $A_i$. Let us choose 
$s_i={1 \over \Delta}$ for $i=1, \ldots, n$ in (1.4.1). Then for the right hand side of (1.4.1), we have 
$$ s_i \prod_{j \in \Gamma_i} (1-s_j)  \ \geq \ {1 \over \Delta} \left(1 - {1 \over \Delta}\right)^{\Delta} \sim {1 \over e \Delta}, $$
and hence our condition (1.2.1) is much more restrictive than (1.4.1).
The constraints (1.4.1) and the lower bound (1.4.2) are, in fact,  sharp \cite{She85} for the intersection $\bigcap_{i=1}^n \overline{A}_i$ to have positive probability,  see also \cite{SS05} and \cite{Jen24}. We, however, are interested in further approximating that probability arbitrarily close.

There has been a lot of work on how to approximate ${\bold P}\left(\bigcap_{i=1}^n \overline{A}_i\right)$ under some conditions, similar to (1.4.1). This includes randomized algorithms, see \cite{Jer24} for a survey of the {\it partial rejection sampling} approach inspired by the algorithmic proof \cite{MT10} of the Lov\'asz Local Lemma, as well as deterministic, see \cite{Moi19}, \cite{JPV22}, \cite{HWY23}, \cite{WY24} and references therein. It appears that they all operate under some additional assumptions, either restricting the intersection pattern of $A_i$  to `extremal' where $A_i \cap A_j =\emptyset$ whenever $A_i$ and $A_j$ share a coordinate \cite{Jer24}, or restricting the structure of the sets $A_i$, such as to those coming from CNF Boolean formulas \cite{Moi19} or atomic clauses with one violation per clause \cite{WY24}, or assuming that each $\Omega_i$ is a finite set with uniform measure and an upper bound on the number of elements in $\Omega_i$ \cite{JPV22}, \cite{HWY23}. 

We note, in particular, that Theorem 1.2 establishes a new regime, where approximate counting of satisfying assignments in the constraint satisfaction problem can be done efficiently, in quasi-polynomial time.

Bencs and Regts \cite{BR24} and then Guo and N \cite{GN25} used an approach somewhat similar to ours to approximate volumes of some combinatorially defined polytopes in deterministic polynomial time (we talk more about the similarities and differences of the approaches in Sections 2.2 and 4.1).

Linial and Nisan \cite{LN90} and then Kahn, Linial and Samorodnitsky \cite{KLS96} considered the problem of approximating ${\bold P} \left(\bigcup_{i=1}^n A_i \right)$ from the probabilities \newline
${\bold P}\left(A_{i_1} \cap \ldots \cap A_{i_k}\right)$ of $k$-wise intersections, without assuming anything about the dependencies of $A_i$. In \cite{LN90}, Linial and Nisan showed in particular 
that for $k = \Omega(\sqrt{n})$, the probability of the union can be approximated from the probabilities of $k$-wise intersections within a relative error of $O\left(e^{-2k/\sqrt{n}}\right)$ and that for $k=O(\sqrt{n})$ the probability of the union cannot be approximated better than within a factor of $O(n/k^2)$. 
Then in \cite{KLS96}, Kahn, Linial and Samorodnitsky proved that for any $k$, the probability of the union can be approximated from the probability of the $k$-wise intersections within an additive error of $e^{-\Omega\left({k^2 \over n \ln n}\right)}$ and that the result is sharp, up to a logarithmic factor in the exponent. Since in our case ${\bold P}\left( \bigcap_{i=1}^n \overline{A}_i\right)$ can be exponentially small in $n$, an attempt to combine the estimates of \cite{LN90} and \cite{KLS96} and the obvious formula (1.1.1) can produce a huge relative error in estimating the probability of the intersection.

The methods of \cite{LN90} and \cite{KLS96} on the one hand, and the method of this paper on the other hand, seem to be quite different. In particular, the final answer is given in different forms: in \cite{LN90} and \cite{KLS96} the value of ${\bold P}\left( \bigcup_{i=1}^n A_i\right)$ is approximated by a linear function in the quantities $\sigma_k$ of (1.2.4),
whereas our Theorem 1.2 approximates $\ln {\bold P} \left( \bigcap_{i=1}^n \overline{A}_i\right)$ by a polynomial in $\sigma_k$. 
Nevertheless, it feels like the approaches must be somehow related. In particular, the wording of one of the results of \cite{KLS96} comes tantalizingly close to the wording of our Theorem 1.2: it is proved in \cite{KLS96} that the number, say $s$,  of satisfying assignments of a Boolean DNF formula in $n$ variables is determined exactly by the set of numbers of satisfying assignments for all conjunctions of $k \leq \log_2 n +1$ terms of the formula, although it is not clear how to compute that $s$ from that set of numbers. 
 
\head 2. The method of the proof of Theorem 1.2 and preliminaries \endhead

\subhead (2.1) The method \endsubhead  Let ${\boldsymbol \Omega}$ be the probability space and let $A_1, \ldots, A_n \subset  {\boldsymbol \Omega}$ be events as in Theorem 1.2. We define a univariate polynomial of a complex variable $z$ by
$$p(z)=\int_{{\boldsymbol \Omega}} \prod_{i=1}^n (1-z [A_i]) \ dx, \tag2.1.1$$
where $[A_i]$ is the indicator of $A_i$ defined by (1.1.3) and $dx$ stands for the integration against the product probability measure in ${\boldsymbol \Omega}$. 
Then 
$$p(1)=\int_{{\boldsymbol \Omega}} \prod_{i=1}^n (1-[A_i]) \ dx = {\bold P}\left( \bigcap_{i=1}^n \overline{A}_i \right).$$
On the other hand, $p(0)=1$ and for the $k$-th derivative of $p$ at $0$ we have
$$\split  p^{(k)}(0)=&(-1)^k k! \sum_{1 \leq i_1 < i_2 < \ldots < i_k \leq n} \int_{{\boldsymbol \Omega}} [A_{i_1}] \cdots [A_{i_k}] \ dx \\= 
&(-1)^k k! \sum_{1 \leq i_1 < i_2 < \ldots < i_k \leq n} {\bold P}\left(A_{i_1} \cap \ldots \cap A_{i_k}\right) \\=
&(-1)^k k! \sigma_k. \endsplit \tag2.1.2$$
Clearly, $\deg p \leq n$. Thus our goal is to approximate $p(1)$ within relative error $\epsilon$ from $p^{(k)}(0)$ with $k =O(\ln (n/\epsilon))$.

\subhead (2.2) Approximating $p(1)$ from $p^{(k)}(0)$ \endsubhead 
It turns out that for our purposes it suffices to show that $p(z) \ne 0$ in some neighborhood of the interval $[0, 1]$ in the complex plane ${\Bbb C}$. More precisely, if a polynomial 
$g(z)$ with $\deg g \leq n$ does not have zeros in a $\delta$-neighborhood of $[0, 1] \subset {\Bbb C}$ for some $0 < \delta < 1$, then for any $0 < \epsilon < 1$, the value of $g(1)$, up to relative error $\epsilon$ (properly defined if $g(1)$ is a non-zero complex number), is determined by $g^{(k)}(0)$ for 
$k \leq K=e^{O(1/\delta)}\ln (n/\epsilon)$.  Given the values $g(0), g'(0), \ldots, g^{(K)}(0)$, the approximation can be computed in $O(K^3)$ time. Furthermore, if $g(0)=1$ then an approximation of $\ln g(1)$ within additive error $\epsilon$ can be written as a polynomial 
$Q$ in $g^{(k)}(0)$ with $k \leq K$ and $\deg Q \leq K$, see Section 2.2 of \cite{Bar16} and also \cite{Bar17}. 

We show that for $\delta=O(1/\Delta)$ and $p(z)$ defined by (2.1.1), we indeed have $p(z) \ne 0$ in the $\delta$-neighborhood of $[0, 1] \subset {\Bbb C}$. 
Assuming that the values of $\sigma_k$ defined by (1.2.4) can be computed in $n^{O(k)}$ time,
the straightforward approach of \cite{Bar16} and \cite{Bar17} allows one to approximate $p(1)={\bold P}\left( \bigcap_{i=1}^n \overline{A}_i\right)$ within relative error 
$\epsilon$ in quasi-polynomial $n^{O\left( \ln (n/\epsilon)\right)}$ time, provided we consider $\Delta$ to be fixed in advance. 

If we assume that ${\bold P}\left(A_{i_1}\cap \ldots \cap A_{i_k}\right)$ can be computed in $2^{O(k)}$ time, there is a possibility of a speed-up to a genuinely polynomial time, using the approach of Patel and Regts \cite{PR17}, see also \cite{LSS19}. Here we very briefly sketch how. With events $A_1, \ldots, A_n$ we associate a hypergraph $H=(V, E)$ with set $V=\{1, \ldots, m\}$ of vertices and set $E \subset 2^V$ of edges, where each event $A_i$ is represented by an edge $J_i \subset \{1, \ldots, m\}$, consisting of the vertices $j$ such that $A_i$ depends on the coordinate $\xi_j$. From the technique of \cite{PR17} and \cite{LSS19} it follows that to compute (2.1.2), it suffices to compute the probabilities when the associated hypergraph is connected (note that for an arbitrary hypergraph, the required probability is the product of the probabilities over its connected components). Then it turns out that 
for a fixed $\Delta$, and for $k=O_{\Delta} \ln (n/\epsilon)$ there are only $\left(n/\epsilon\right)^{O_{\Delta}(1)}$ connected subhypergraphs of $H$ with $k$ edges.
Very recently, a different way to accelerate computations was described in \cite{BR25}.

We note that a particular version of the polynomial (2.1.1) was used earlier by Bencs and Regts \cite{BR24} to approximate the volume of a truncated independence polytope of a graph and then by Guo and N \cite{GN25} to approximate the volume of a truncated fractional matching polytope.

Thus our goal is to prove that the polynomial $p(z)$ of (2.1.1) does not have zeros in a neighborhood of the interval $[0, 1] \subset {\Bbb C}$.
We let 
$$U=\Bigl\{z \in {\Bbb C}: \ \operatorname{dist}(z, [0,1]) \leq \delta \Bigr\} \quad \text{for} \quad \delta={1 \over 6\Delta},$$
where $\operatorname{dist}$ is the Euclidean distance in the complex plane ${\Bbb C}={\Bbb R}^2$.

Let us fix a $z \in U$ and let us define $f_i: {\boldsymbol \Omega} \longrightarrow {\Bbb C}$ by 
$$f_i=1-z[A_i] \quad \quad \text{for} \quad i=1, \ldots, n. $$ Then 
$$|f_i(x)| \ \leq \ 1+ {1 \over 6\Delta} \quad \text{for all} \quad x \in {\boldsymbol \Omega} \quad \text{and} \quad i=1, \ldots, n $$
and 
$${\bold P}\left(f_i(x) \ne 1 \right) \ < \  (3 \Delta)^{-3\mu_i}\quad \text{for} \quad i=1, \ldots, n. $$
Hence we will be proving the following reformulation of Theorem 1.2.

\proclaim{(2.3) Theorem} Let $\Omega_1, \ldots, \Omega_m$ be probability spaces and let ${\boldsymbol \Omega}=\Omega_1 \times \cdots \times \Omega_m$ be their product. Let $f_1, \ldots, f_n: {\boldsymbol \Omega} \longrightarrow {\Bbb C}$ be random variables, such that each $f_i$ depends on not more than $r_i$ coordinates 
and shares a coordinate with at most $\Delta_i$ of other functions $f_j$. Let 
$$\Delta=\max\left\{5, \ \Delta_i: \ i=1, \ldots, n \right\}$$
 and let 
$$\mu_i =\min\left\{ r_i,\ \Delta_i+1\right\} \quad \text{for} \quad i=1, \ldots, n.$$
Suppose further that
$$\left| f_i(x)\right| \ \leq \ 1 + {1 \over {6 \Delta}} \quad \text{for all} \quad x \in {\boldsymbol \Omega} \quad \text{and} \quad i=1, \ldots, n \tag2.3.1$$
and that
$${\bold P}\left(f_i \ne 1 \right) \ < \ (3 \Delta)^{-3 \mu_i} \quad \text{for} \quad i=1, \ldots, n. \tag2.3.2$$
Then
$$\int_{{\boldsymbol \Omega}} \prod_{i=1}^n f_i(x) \ d x \ne 0.$$
\endproclaim

We will prove Theorem 2.3  by proving a stronger statement by induction on the number $n$ of functions.

We will use some simple geometric facts.
\proclaim{(2.4) Lemma} Let $\Omega$ be a probability space and let $f: \Omega \longrightarrow {\Bbb C}$ be an integrable random variable. Suppose that for every $\omega \in \Omega$ we have $f(\omega) \ne 0$ and that for any $\omega_1, \omega_2 \in \Omega$ the angle between $f(\omega_1) \ne 0$ and $f(\omega_2) \ne 0$, 
considered as vectors in ${\Bbb R}^2 = {\Bbb C}$, does not exceed $\theta$ for some $0 \leq \theta < 2\pi/3$. Then 
$$\left| \int_{\Omega} f(x) \ dx \right|  \ \geq \ \left( \cos {\theta \over 2}\right) \int_{\Omega} |f(x)| \ dx.$$
\endproclaim
\demo{Proof} See Lemma 3.3 of \cite{BR19} and also Lemma 3.6.3 of \cite{Bar16}.
{\hfill \hfill \hfill} \qed
\enddemo

\proclaim{(2.5) Lemma} Let $v \in {\Bbb C} \setminus \{0\}$ be a non-zero complex number and let $u \in {\Bbb C}$ be another complex number such that 
$$|u| \ \leq \ (\sin \theta) |v| \quad \text{for some} \quad 0 \leq \theta < \pi/2.$$
Then the angle between $v$ and $v+u$, considered as vectors in ${\Bbb R}^2 ={\Bbb C}$, does not exceed $\theta$. 
\endproclaim
\demo{Proof} See Lemma 3.6.4 of \cite{Bar16}.
{\hfill \hfill \hfill} \qed
\enddemo

We will also use the classical Markov inequality in the following situation.

\proclaim{(2.6) Lemma} Let $\Omega_1$ and $\Omega_2$ be probability spaces with respective probability measures ${\bold P}_1$ and ${\bold P}_2$, let $\Omega=\Omega_1 \times \Omega_2$ be their product with probability measure ${\bold P}={\bold P}_1 \times {\bold P}_2$, and let 
$A \subset \Omega$ be an event. For $\omega_1 \in \Omega_1$ we define an event $A_{\omega_1} \subset \Omega_2$ by 
$$A_{\omega_1}=\left\{ \omega_2 \in \Omega_2: \quad (\omega_1, \omega_2) \in A \right\}.$$
Then for any $\tau \geq 1$, we have 
$${\bold P}_1 \Bigl( \omega_1: \ {\bold P}_2 \left( A_{\omega_1}\right) \ \geq \ \tau {\bold P}(A) \Bigr) \ \leq \ \tau^{-1}.$$
\endproclaim
\demo{Proof} We define a random variable $f: \Omega_1 \longrightarrow [0, 1]$ by 
$$f(\omega_1) = {\bold P_2} \left(A_{\omega_1}\right) \quad \text{for} \quad \omega_1 \in \Omega_1.$$
Then 
$$\int_{\Omega_1} f(\omega_1) \ d \omega_1 = {\bold P}(A),$$
where $d\omega_1$ stands for the integration against the probability measure ${\bold P_1}$.
Since $f$ is non-negative, the result follows.
{\hfill \hfill \hfill} \qed
\enddemo

\head 3. Proof of Theorem 2.3 \endhead

Our goal is to prove Theorem 2.3. We do that by proving a stronger result by induction on the number $n$ of functions $f_i$.

\subhead (3.1) The induction hypothesis \endsubhead Let
$$\theta={1 \over  \Delta^2} \quad \text{and} \quad \beta = 1+{1 \over \Delta}.$$

The induction hypothesis consists of the following two statements. 
\bigskip
$\bullet$ {\sl Statement I:} Let $f_i: {\boldsymbol \Omega} \longrightarrow {\Bbb C}$, $i=1, \ldots, n$, be random variables as in Theorem 2.3.
Then
$$\int_{{\boldsymbol \Omega}} \prod_{i=1}^n f_i(x) \ dx \ne 0, \quad \int_{{\boldsymbol \Omega}} \prod_{i=1}^{n-1} f_i(x) \ dx \ne 0 \tag3.1.1$$
and the angle between two non-zero complex numbers (3.1.1), considered as vectors in ${\Bbb R}^2 = {\Bbb C}$, does not exceed $\theta$.
For $n=1$, we agree that the second number is 1, following the general convention that the product of the numbers from an empty set is 1, while the sum of the numbers from an empty set is 0.
\bigskip
$\bullet$ {\sl Statement II:} Let $f_i: {\boldsymbol \Omega} \longrightarrow {\Bbb C}$, $i=1, \ldots, n$, be random variables as in Statement I.
For $k \leq n$, let $\widehat{f}_{n-k+1}, \ldots, \widehat{f}_n: {\boldsymbol \Omega} \longrightarrow {\Bbb C}$ be random variables such 
that 
$$\left| \widehat{f}_i(x)\right| \ \leq \ 1+ {1 \over 6\Delta} \quad \text{for all} \quad x \in {\boldsymbol \Omega} \quad \text{and} \quad i=n-k+1, \ldots, n, \tag3.1.2$$
and such that each $\widehat{f}_i$ depends on a subset of the coordinates that  $f_i$ depends on (note, however, that $\widehat{f}_i$ do not have to satisfy (2.3.2)). Then 
$$\left| \int_{{\boldsymbol \Omega}} \left(\prod_{i=1}^{n-k} f_i(x) \prod_{i=n-k+1}^n \widehat{f}_i(x) \right)\ dx \right| \ \leq \ \beta^k \left| \int_{{\boldsymbol \Omega}}
\prod_{i=1}^n f_i(x) \ dx \right|. $$
In words: if we replace $k$ of the functions $f_i$ by some functions $\widehat{f}_i$, depending on the same or smaller set of coordinates, and still satisfying (2.3.1) but not necessarily (2.3.2), then the absolute value of the integral of the product can increase by a factor $\beta^k$ at most.
\bigskip
\remark{(3.2) Remarks}  Suppose that Statement I holds, and let $\widehat{f}_n: {\boldsymbol \Omega} \longrightarrow {\Bbb C}$ be yet another random variable that depends on a subset of coordinates that $f_n$ depends on and such that 
$$\left| \widehat{f}_n(x)\right| \ \leq \ 1 + {1 \over 6 \Delta} \quad \text{and} \quad {\bold P}\left( \widehat{f}_n(x) \ne 1 \right) \ < \ (3 \Delta)^{-3\mu_n}.$$
Then the angle between two complex numbers
$$\int_{{\boldsymbol \Omega}} \prod_{i=1}^n f_i(x) \ dx \ne 0 \quad \text{and} \quad \int_{{\boldsymbol \Omega}} \widehat{f}_n(x) \prod_{i=1}^{n-1} f_i(x) \ dx \ne 0 \tag3.2.1$$
does not exceed $\displaystyle 2 \theta={2 \over \Delta^2}$. Indeed, both numbers (3.2.1) form an angle of at most $\theta$ with 
$\displaystyle \int_{{\boldsymbol \Omega}} \prod_{i=1}^{n-1} f_i(x) \ dx \ne 0$.
Iterating, we conclude that if $\widehat{f}_{n-k+1}, \ldots, \widehat{f}_n: {\boldsymbol \Omega} \longrightarrow {\Bbb C}$ are random variables satisfying (3.1.2) and also such that 
$${\bold P}\left(\widehat{f}_i \ne 1 \right) \ < \ (3\Delta)^{-3\mu_i} \quad \text{for} \quad i=n-k+1, \ldots, n$$
then the angle between two complex numbers
$$\int_{{\boldsymbol \Omega}} \prod_{i=1}^n f_i(x) \ dx \ne 0 \quad \text{and} \quad \int_{{\boldsymbol \Omega}} \left( \prod_{i=1}^{n-k} f_i(x) \prod_{i=n-k+1}^n \widehat{f}_i(x) \right) \ dx \ne 0$$
does not exceed $2 k \theta = 2k/\Delta^2$. 

In words: if we replace $k$ of the random variables $f_i$ satisfying (2.3.1) and (2.3.2) by some other random variables on the same or smaller sets of coordinates satisfying the same conditions, then the value of the integral rotates by at most an angle of $2k/\Delta^2$.

Finally, for the sake of transparency,  in the proof below we replace condition (2.3.2) by 
$${\bold P}(f_i \ne 1) \ < \ (\gamma \Delta)^{-3\mu_i}, \tag3.2.2$$
where $\gamma > 0$ is some absolute constant. In the course of the proof it then becomes clear that the induction can be carried through for a sufficiently large $\Delta$ and $\gamma$, and numerical computations show that we can indeed choose $\Delta \geq 5$ and $\gamma=3$, as claimed.
\endremark

\subhead (3.3) Base $n=1$ \endsubhead In this case, we have just one function $f_1$ in some  $r \leq r_1$ coordinates, satisfying (2.3.1) and (3.2.2). 

If $r=0$ then $\mu_1=0$. Since $r=0$, the function $f_1$ is a constant and by (3.2.2), we must have $f_1 \equiv 1$, so that 
$$\int_{{\boldsymbol \Omega}} f_1(x) \ dx =1.$$
Statement I then holds tautologically.  By (3.1.2), 
$$\left| \int_{{\boldsymbol \Omega}} \widehat{f}_1(x)  \ d x \right| \ \leq \ 1+ {1 \over 6\Delta} \ \leq \ 
\beta \left| \int_{{\boldsymbol \Omega}} f_1(x) \ dx \right|,$$
and hence Statement II also holds.

Suppose now that $r \geq 1$, so that $\mu_1=1$. Let 
$$B=\Bigl\{x \in {\boldsymbol \Omega}: \quad f_1(x)=1 \Bigr\}, \quad \text{so that} \quad {\bold P}(\overline{B}) \ < \ (\gamma \Delta)^{-3}.$$
From (2.3.1) and (3.2.2), 
$$\aligned \left|1-\int_{{\boldsymbol \Omega}} f_1(x) \ dx \right| = &\left| \int_{\overline{B}} (1- f_1(x)) \ dx \right|
 \ \leq \ \left(2 + {1 \over 6 \Delta}\right) ( \gamma \Delta)^{-3}. \endaligned
\tag3.3.1$$ 
It is clear now that for $\gamma \geq 1$ sufficiently large, we have 
$$ \left(2 + {1 \over 6 \Delta}\right) ( \gamma \Delta)^{-3} \ \leq \ \sin {1 \over \Delta^2} =\sin \theta, \tag3.3.2$$
so that by Lemma 2.5 the angle between complex numbers $1$ and  $\int_{{\boldsymbol \Omega}} f_1(x) \ dx  \ne 0$  does not exceed $\theta$ and Statement I holds. 
Numerical computations show that (3.3.2) indeed holds for $\Delta \geq 3$ and $\gamma \geq 1$.

To check Statement II, from  (3.3.1), we have 
$$\left| \int_{{\boldsymbol \Omega}} f_1(x)\ dx \right| \ \geq \ 1 - \left(2 + {1 \over 6\Delta}\right) (\gamma \Delta)^{-3}. \tag3.3.3$$
From (3.1.2), we have 
$$ \left| \int_{{\boldsymbol \Omega}} \widehat{f}_1(x)  \ d x \right| \ \leq \ 1+ {1 \over 6\Delta}. \tag3.3.4$$
It is now clear that for $\gamma > 1$ sufficiently large, we have 
$$1+ {1 \over 6\Delta} \ \leq \ \left(1+{1 \over \Delta}\right)  \left(1 - \left(2 + {1 \over 6\Delta}\right) (\gamma \Delta)^{-3}\right), \tag3.3.5$$
and so by (3.3.3) and (3.3.4), 
$$\left| \int_{{\boldsymbol \Omega}} \widehat{f}_1(x) \ dx \right| \ \leq \ \beta \left| \int_{{\boldsymbol \Omega}} f_1(x) \ dx \right|,$$
and  Statement II holds too.

Numerical computations show that (3.3.5) holds provided $\Delta \geq 3$ and $\gamma \geq 1$.

\subhead (3.4) Induction step $n-1 \Longrightarrow n$, trivial case \endsubhead Thus we assume that $n \geq 2$. Suppose that 
$f_n$ depends on the coordinates $\xi_j: \ j \in J$, so $|J| =r \leq r_n$. If $r=0$ then $f_n$ is a constant, and by (3.2.2), we have $f_n \equiv 1$. 
Then 
$$\prod_{i=1}^n f_i(x)= \prod_{i=1}^{n-1} f_i(x)$$
and Statement I holds trivially.

Since $\widehat{f}_n$ depends on a subset of the set of coordinates that $f_n$ depends on, the function $\widehat{f}_n$ is also a constant. Then by (3.1.2)
$$\split &\left| \int_{{\boldsymbol \Omega}} \left( \prod_{i=1}^{n-k} f_i(x) \prod_{i=n-k+1}^n \widehat{f}_i(x)\right) \ d x\right| \\ \leq \ 
&\left(1+ {1 \over 6\Delta}\right)  \left| \int_{{\boldsymbol \Omega}} \left( \prod_{i=1}^{n-k} f_i(x) \prod_{i=n-k+1}^{n-1} \widehat{f}_i(x) \right) \  d x\right| \\
\leq \ & \beta^k \left| \int_{{\boldsymbol \Omega}} \prod_{i=1}^n f_i(x)\ dx \right|,
\endsplit$$
where the last inequality follows by the induction hypothesis. Thus Statement II also holds.

\subhead (3.5) Induction step $n-1 \Longrightarrow n$, preparations for the general case \endsubhead 
Without loss of generality, we assume now that the set $J$ of coordinates that $f_n$ depends on is non-empty, so $f_n$ depends on some $r$ coordinates,
$1 \leq r \leq r_n$.
We form two product probability spaces 
$${\boldsymbol \Omega}_1 = \times_{j \in J} \Omega_j \quad \text{and} \quad {\boldsymbol \Omega}_2 = \times_{j \notin J} \Omega_j, $$
so that 
$${\boldsymbol \Omega} = {\boldsymbol \Omega}_1 \times {\boldsymbol \Omega}_2.$$
We then write $x \in {\boldsymbol \Omega}$ as $x=(u, y)$, where $u \in {\boldsymbol \Omega_1}$ and $y \in {\boldsymbol \Omega_2}$. 
We denote the probability measure in ${\boldsymbol \Omega}$ by ${\bold P}$, in ${\boldsymbol \Omega_1}$ by ${\bold P}_1$ and in ${\boldsymbol \Omega}_2$ by 
${\bold P}_2$. We use the notation $dx$, $du$ and $dy$ to denote integration against ${\bold P}$, ${\bold P}_1$ and ${\bold P}_2$ respectively.

Let $I \subset \{1, \ldots, n-1\}$ be the set of indices $i$ such that $f_i(x)$ shares a coordinate with $f_n(x)$. 
Then 
$$|I| \leq \Delta. \tag3.5.1$$
For $i \in I$ and $u \in {\boldsymbol \Omega_1}$, we define a random variable $h_i(\cdot |u): {\boldsymbol \Omega_2} \longrightarrow {\Bbb C}$ by 
$$h_i(y| u) = f_i(u, y) \quad \text{for} \quad u \in {\boldsymbol \Omega_1}, \quad y \in {\boldsymbol \Omega}_2 \quad \text{and} \quad i \in I. $$

Similarly, for $i \in I$ such that $i > n-k$ and random variables $\widehat{f}_i$ as in Statement II of Section 3.1, we define a random variable $\widehat{h}_i(\cdot| u): {\boldsymbol \Omega_2} \longrightarrow {\Bbb C}$ by
$$\widehat{h}_i(y| u) = \widehat{f}_i(u, y) \quad \text{for} \quad u \in {\boldsymbol \Omega_1}, \quad y \in {\boldsymbol \Omega}_2 \quad \text{and} \quad i \in I. $$
We note that each $h_i(y|u)$ depends on at most $r_i-1$ coordinates and shares a coordinate with at most $\Delta_i-1$ of functions 
$h_j(y|u)$ for $j \in I \setminus \{i\}$ and $f_j(y)$ for $j \notin I \cup \{n\}$.

We further define functions $\phi, \widehat{\phi}: {\boldsymbol \Omega_1} \longrightarrow {\Bbb C}$ by 
$$\aligned &\phi(u) = \int_{{\boldsymbol \Omega_2}} \left( \prod_{i \in I} h_i(y|u) \prod\Sb i \notin I \\ i < n \endSb f_i(y) \right)\ dy
 \quad \text{and} \\
&\widehat{\phi}(u) =\int_{{\boldsymbol \Omega_2}} \left(\prod\Sb i \in I \\ i \leq n-k \endSb h_i(y|u)
 \prod\Sb i \notin I \\ i \leq n-k \endSb f_i(y) \prod\Sb i \in I \\ i > n-k \endSb \widehat{h}_i(y|u)  \prod\Sb i \notin I \\ n-k < i < n \endSb \widehat{f}_i(y) \right)\ dy.
\endaligned
\tag 3.5.2$$
Then by Fubini's Theorem
$$\split & \int_{{\boldsymbol \Omega}} \prod_{i=1}^n f_i(x) \ dx = \int_{{\boldsymbol \Omega_1}} f_n(u) \phi(u) \ du \\
&\int_{{\boldsymbol \Omega}} \prod_{i=1}^{n-1} f_i(x) \ dx = \int_{{\boldsymbol \Omega_1}} \phi(u) \ du  \quad \text{and} \\
&\int_{{\boldsymbol \Omega}} \left( \prod_{i=1}^{n-k} f_i(x) \prod_{i=n-k+1}^n \widehat{f}_i(x)\right) \ dx = \int_{\boldsymbol \Omega_1} \widehat{f}_n(u) \widehat{\phi}(u) \ du.
\endsplit \tag3.5.3$$
Applying Lemma 2.6 with $\tau=(\gamma \Delta)^3$, from (3.2.2),  we conclude that 
$${\bold P_1} \Bigl(u: \ {\bold P_2} (h_i(\cdot | u) \ne 1 ) \ \geq \ (\gamma \Delta)^{-3(\mu_i-1)} \Bigr) \ \leq \ (\gamma \Delta)^{-3} \quad \text{for} \quad i \in I. \tag3.5.4$$ 
We define an event $A \subset {\boldsymbol \Omega_1}$ by
$$A=\Bigl\{u: \ {\bold P}_2(h_i(\cdot | u) \ne 1 ) \ < \ (\gamma \Delta)^{-3(\mu_i-1)} \quad \text{for all} \quad i \in I \Bigr\}. \tag3.5.5$$
Then by (3.5.1) and (3.5.4)
$${\bold P}_1(\overline{A}) \ \leq \ {1 \over \gamma^3 \Delta^2}. \tag3.5.6$$
We define an event $B \subset {\boldsymbol \Omega_1}$ by 
$$B=\Bigl\{u:\ f_n(u) =1 \Bigr\}. \tag3.5.7$$
so that by (3.2.2),
$${\bold P}_1(\overline{B}) \ < \ (\gamma \Delta)^{-3\mu_n} \ \leq \ (\gamma \Delta)^{-3}.\tag3.5.8$$

Since each function $h_i(\cdot | u)$ depends on at most $r_i-1$ coordinates and shares a coordinate with at most $\Delta_i-1$ functions $h_j(\cdot | u)$ for $j \in I \setminus \{i\}$ and 
$f_j(y)$ for $j \in \{1, \ldots, n-1\}\setminus I$, by (3.5.5) for all $u \in A$ the $n-1$ functions $h_i(\cdot | u)$ for $i \in I$ and $f_i$ for $i \in \{1, \ldots, n-1\} \setminus I$ satisfy Statement I of the induction hypothesis. Let us choose any two $u_1, u_2 \in A$.  When we switch from $u_1$ to $u_2$ in the integral (3.5.2) for $\phi$, we change the $|I|$ functions $h_i(\cdot | u_1)$ to $h_i(\cdot | u_2)$, and hence applying Statement I of the induction hypothesis and Remark 3.2, we conclude that $\phi(u_1) \ne 0$, $\phi(u_2) \ne 0$ and the angle between these two non-zero complex numbers is at most $2 |I| \theta \ \leq \ 2/\Delta$  by (3.5.1). Hence by Lemma 2.4,
$$\left| \int_{A \cap B} \phi(u) \ du \right| \ \geq \ \left(\cos {1 \over \Delta}\right) \int_{A \cap B} |\phi(u)| \ du \ > \ 0. \tag3.5.9$$

Let us now pick any $u_1 \in A$ and any $u_2 \in {\boldsymbol \Omega_1}$. If in the formula (3.5.2) for $\phi(u)$, we switch $u_1$ to $u_2$, we switch at most $|I| \leq \Delta$ of functions 
$h_i(\cdot |u_1)$ to $h_i(\cdot|u_2)$. Applying Statement II of the induction hypothesis, we conclude that 
$$\aligned |\phi(u_2)| \ \leq \ &\beta^{|I|}  |\phi(u_1)| \ \leq \ \left(1 + {1 \over \Delta}\right)^{\Delta}  |\phi(u_1)| \ \leq \ 3 |\phi(u_1)| \\
&\text{for all} \quad u_1 \in A, \ u_2 \in {\boldsymbol \Omega_1}. \endaligned \tag3.5.10$$
Similarly, to pass from $\phi(u)$ to $\widehat{\phi}(u)$ for $u \in A$, we change at most $k-1$ in total functions $h_i(\cdot| u)$ to $\widehat{h}_i(\cdot | u)$ and
$f_i(y)$ to $\widehat{f}_i(y)$ in (3.5.2), and hence applying Statement II of the induction hypothesis, we obtain
$$|\widehat{\phi}(u)| \ \leq \ \beta^{k-1} |\phi(u)| \quad \text{for all} \quad u \in A. \tag3.5.11$$
On the other hand, in $\widehat{\phi}(u_2)$ for $u_2 \in {\boldsymbol \Omega_1}$ compared to $\phi(u_1)$ for $u_1 \in A$, we replace 
additionally up to $|I| \leq \Delta$ functions $h_i(\cdot |u_1)$ by $h_i(\cdot | u_2)$, and so from Statement II of the induction hypothesis, we get 
$$\split |\widehat{\phi}(u_2)| \ \leq \ &\beta^{k-1} \beta^{|I|}  |\phi(u_1)| \ \leq \ 3\beta^{k-1} |\phi(u_1)| \\
 &\text{for all} \quad u_1 \in A, \ u_2 \in {\boldsymbol \Omega_1}. \endsplit \tag3.5.12$$

\subhead (3.6) Induction step $n-1 \Longrightarrow n$, Statement I \endsubhead We compare the first two integrals in (3.5.3) and argue that the bulk of both integrals come from integrating over the event $A \cap B$, for $A$ defined by (3.5.5) and $B$ defined by (3.5.7), where the integrals obviously coincide.
Indeed, by (3.5.7) we have 
$$\int_{A \cap B} f_n(u) \phi(u) \ du= \int_{A \cap B} \phi(u)\ du.$$ 
From (3.5.3), 
$$\left| \int_{\boldsymbol{\Omega}} \prod_{i=1}^n f_i(x) \ d x - \int_{A \cap B} \phi(u) \ du \right| \ \leq \  \int_{\overline{A} \cup \overline{B}} | f_n(u) \phi(u)| \ du 
 \tag 3.6.1$$
and similarly,
$$\left| \int_{\boldsymbol{\Omega}} \prod_{i=1}^{n-1} f_i(x) \ d x - \int_{A \cap B} \phi(u) \ du \right| \ \leq \  \int_{\overline{A} \cup \overline{B}} |  \phi(u)| \ du. \tag 3.6.2$$
 Using (3.1.2), (3.5.10), (3.5.6), (3.5.8)  and  (3.5.9) in that order, we conclude that 
 $$\aligned &\int_{\overline{A} \cup \overline{B}} |f_n(u) \phi(u)| \ d u \ \leq \ \left(1 + {1 \over 6 \Delta}\right) \int_{\overline{A} \cup \overline{B}} |\phi(u)| \ du \\ &\qquad \leq \ 
\left(1+ {1 \over 6 \Delta}\right) { 3{\bold P}_1(\overline{A} \cup \overline{B}) \over {\bold P}_1(A \cap B)} \int_{A \cap B} |\phi(u)| \ du \\ &\qquad \leq
 \ \left(1 + {1 \over 6\Delta}\right) \left({6 \over \gamma^3 \Delta^2- 2 }\right) \int_{A \cap B} |\phi(u)| \ du \\
 &\qquad \leq \ \left(1 + {1 \over 6\Delta}\right) \left({6 \over \gamma^3 \Delta^2 -2}\right)\left( {1 \over \cos(1/\Delta)}\right) \left| \int_{A \cap B} \phi(u) \ du \right|. \endaligned \tag3.6.3$$
 and, similarly,
$$ \int_{\overline{A} \cup \overline{B}} |\phi(u)| \ d u \ \leq \ \left({6 \over \gamma^3 \Delta^2 -2}\right)\left( {1 \over \cos(1/\Delta)}\right) \left| \int_{A \cap B} \phi(u) \ du \right|. \tag3.6.4$$
In view of (3.6.1) -- (3.6.4) and Lemma 2.5, to make sure that 
$$\int_{{\boldsymbol \Omega}} \prod_{i=1}^n f_i(x) \ dx \ne 0, \quad \int_{{\boldsymbol \Omega}} \prod_{i=1}^{n-1} f_i(x) \ dx \ne 0$$
and that the angles between each of these two integrals and 
$$ \int_{A \cap B} \phi(u) \ du \ \ne 0,$$
cf. (3.5.9), do not exceed $\theta/2 = {1 \over 2 \Delta^2}$ each, we should make sure that 
$$ \left(1 + {1 \over 6\Delta}\right) \left({6 \over \gamma^3 \Delta^2 -2}\right)\left( {1 \over \cos(1/\Delta)}\right) \ \leq \ \sin {1 \over 2 \Delta^2}. \tag3.6.5$$
Clearly, (3.6.5) holds for a sufficiently large $\gamma > 1$. Computations show that (3.6.5) holds
when $\Delta \geq 3$ and $\gamma \geq 3$, so Statement I holds.
 
 \subhead (3.7) Induction step $n-1 \Longrightarrow n$, Statement II \endsubhead Here we compare the first and the last integrals in (3.5.3), and once again argue that the main contribution to both come from the integral over the event $A \cap B$, for $A$ defined by (3.5.5) and $B$ defined by (3.5.7). From (3.5.3) and (3.1.2), we have 
 $$\aligned &\left| \int_{{\boldsymbol \Omega}} \left( \prod_{i=1}^{n-k} f_i(x) \prod_{i=n-k+1}^n \widehat{f}_i(x) \right) \ dx\right| = \left| \int_{{\boldsymbol \Omega_1}} \widehat{f}_n(u) \widehat{\phi}(u) \ du \right| \\
 &\qquad \leq \ \left(1+ {1 \over 6\Delta}\right) \int_{{\boldsymbol \Omega_1}} |\widehat{\phi}(u)| \ du. \endaligned \tag3.7.1$$
 As before, the last integral we split into two cases: over $A \cap B$ and over $\overline{A} \cup \overline{B}$.
 From (3.5.11), we have 
 $$\int_{A \cap B} |\widehat{\phi}(u)| \ d u \ \leq \ \beta^{k-1} \int_{A \cap B} |\phi(u)| \ du \tag 3.7.2$$
 and from (3.5.12), we have 
 $$\int_{\overline{A} \cup \overline{B}} |\widehat{\phi}(u)| \ du \ \leq \ 3 \beta^{k-1} {{\bold P}_1(\overline{A} \cup \overline{B}) \over  {\bold P}_1(A \cap B)} 
 \int_{A \cap B} |\phi(u)| \ du. \tag3.7.3$$
 Hence from  (3.7.1)--(3.7.3) as well as (3.5.6), (3.5.8) and (3.5.9), we get 
 $$\aligned &\left| \int_{{\boldsymbol \Omega_1}} \widehat{f}_n(u) \widehat{\phi}(u) \ du \right| 
 \ \leq \ \beta^{k-1} \left(1 + {1 \over 6 \Delta}\right) \\ &\qquad \times \left(1 + 3  {{\bold P}_1(\overline{A} \cup \overline{B}) \over  {\bold P}_1(A \cap B)} \right) \int_{A \cap B} |\phi(u)| \ du \\ 
&\qquad \leq \ \beta^{k-1}\left(1+ {1 \over 6 \Delta}\right) \left(1+{6 \over \gamma^3 \Delta^2 -2}\right)\left( {1 \over \cos(1/\Delta)}\right) \left| \int_{A \cap B} \phi(u) \ du \right|.\endaligned $$
It is clear now that for sufficiently large $\Delta \geq 1$ and $\gamma > 1$, we have 
$$\left(1+ {1 \over 6 \Delta}\right) \left(1+{6 \over \gamma^3 \Delta^2 -2}\right)\left( {1 \over \cos(1/\Delta)}\right) \ \leq \ 1 + {1 \over 3 \Delta} \tag3.7.4$$
and computations show that (3.7.4) holds for $\Delta \geq 5$ and $\gamma \geq 3$. Hence
$$\split &\left| \int_{{\boldsymbol \Omega_1}} \widehat{f}_n(u) \widehat{\phi}(u) \ du \right| \ \leq \ \beta^{k-1} \left(1 + {1 \over 3 \Delta}\right) \left| \int_{A \cap B} \phi(u) \ du\right| \\ &\qquad \text{provided} \quad \Delta \geq 5 \quad \text{and} \quad \gamma \geq 3. \endsplit \tag3.7.5$$
On the other hand, from (3.5.3) and (3.6.3),
$$\aligned &\left| \int_{\Omega} \prod_{i=1}^n f_i(x) \ dx\right| = \left| \int_{{\boldsymbol \Omega_1}} f_n(u) \phi(u) \ du \right| \\ &\qquad\geq \ \left| \int_{A \cap B} \phi(u) \ du \right| - \left| \int_{\overline{A} \cup \overline{B}} f_n(u) \phi(u) \ du \right|\\ 
&\qquad \geq \ \left(1- \left(1 + {1 \over 6\Delta}\right) \left({6 \over \gamma^3 \Delta^2 -2}\right)\left( {1 \over \cos(1/\Delta)}\right)\right) \left| \int_{A \cap B} \phi(u) \ du \right|. \endaligned $$
It is now clear that if $\Delta \geq 1$ and $\gamma \geq 1$ are sufficiently large, then 
$$\left(1+{1 \over 6\Delta}\right) \left({6 \over \gamma^3 \Delta^2 -2}\right)\left( {1 \over \cos(1/\Delta)}\right) \ \leq \ {1 \over 15 \Delta}. \tag3.7.6$$
Computations show that (3.7.6) indeed holds for $\Delta \geq 5$ and $\gamma \geq 3$, so we have
$$\split &\left| \int_{{\boldsymbol \Omega_1}} f_n(u) \phi(u) \ du \right| \ \geq \ \left(1 - {1 \over 15 \Delta}\right) \left| \int_{A\cap B} \phi(u) \ du \right| \\
&\quad \text{provided} \quad \Delta \geq 5 \quad \text{and} \quad \gamma \geq 3.\endsplit \tag3.7.7$$
Since for $\Delta \geq 5$, we have 
$$\left( 1+ {1 \over 3 \Delta}\right) \left(1- {1 \over 15 \Delta}\right)^{-1} \ < \ 1+ {1 \over \Delta} =\beta,$$
combining (3.5.3), (3.7.5) and (3.7.7), we conclude that Statement II holds.
{\hfill \hfill \hfill} \qed 

\head 4. Concluding remarks \endhead

\subhead (4.1) Comparison with the cluster expansion method \endsubhead One can prove a version of Theorem 1.2 using the method of {\it cluster expansion}, see \cite{Jen24} for a survey. This method was also used by Bencs and Regts \cite{BR24} in the particular case of approximating the volume of a truncated independence polytope of a graph and then by Guo and N \cite{GN25} to approximate the volume of a truncated  matching polytope of a hypergraph. We briefly sketch how it could work in the general framework of Theorem 1.2. 

Given events $A_1, \ldots, A_n$ in as in Theorem 1.2, we still consider the polynomial
$$\split p(z)=&\int_{{\boldsymbol \Omega}} \prod_{i=1}^n (1-z [A_i]) \ dx\\= &\sum_{k=0}^{n} (-1)^k z^k \sum_{1 \leq i_1 < \ldots < i_k \leq n} {\bold P}\left(A_{i_1} \cap \ldots \cap A_{i_k}\right) \endsplit
\tag4.1.1$$
and aim to prove that $p(z) \ne 0$ in a neighborhood of $[0, 1] \subset {\Bbb C}$, see Sections 2.1 and 2.2.
We construct a graph $G=(V, E)$ with set $V=\{1, \ldots, n\}$ of vertices, and vertices $i$ and $j$ spanning an edge whenever  $A_i$ and $A_j$ share a coordinate. The expansion (4.1.1) allows one to apply the {\it polymer model} here: for each collection $1 \leq i_1 < \ldots < i_k \leq n$, we consider the subgraph $H=H_{i_1, \ldots, i_k}$ of $G$ induced by $i_1, \ldots, i_k$ and define its (complex) weight by
$$w\left(H_{i_1, \ldots, i_k}\right) = (-1)^k z^k {\bold P}\left(A_{i_1} \cap \ldots \cap A_{i_k}\right).$$
Then (4.1.1) is the sum of $w(H)$ over all induced subgraphs $H$ of $G$, including the empty subgraph with weight 1.

We say that two connected induced subgraphs $H_1$ and $H_2$ of $G$ are {\it compatible}, denoted $H_1 \sim H_2$,  if they are vertex-disjoint and there are no vertices $v_1$ of $H_1$ and $v_2$ of $H_2$ spanning an edge of $G$, and {\it incompatible}, denoted $H_1 \not\sim H_2$, otherwise.
The crucial property that makes the polymer model work is that if the induced subgraph $H$ is a union of pairwise compatible connected subgraphs, say $H^1, \ldots, H^s$, then 
$$w(H) = w(H^1) \cdots w(H^s).$$
The Koteck\'y - Preiss condition \cite{KP86},  see also \cite{Jen24} for a survey, then provides a sufficient condition for $p(z)$ to be non-zero: if one can assign non-negative real weights $\psi(H) \geq 0$ to the induced subgraphs $H$ such that for every connected induced subgraph $F$ of $G$, we have 
$$\sum_{H:\ H \not\sim F} |w(H)| e^{\psi(H)} \ \leq \ \psi(F), \tag 4.1.2$$
where the sum is taken over all connected induced subgraphs of $G$ that are incompatible with $F$, including $F$ itself. One can bound $|w(H)|$ from above by arguing that the graph $H$ with $|H|$ vertices contains a subset of size at least $|H|/(\Delta+1)$, corresponding to mutually independent events.
The condition (4.1.2) allows one to establish that $p(z) \ne 0$ 
in a disc
$$\bigl\{ z \in {\Bbb C}: \quad |z| < \rho \bigr\} \tag4.1.3$$
of some radius $\rho > 0$ (one should choose $\psi(H)$ to be proportional to the number $|H|$ of vertices in $H$).
The bound one gets from this approach is 
$${\bold P}(A_i) < \Delta^{-\gamma \Delta},\tag4.1.4$$ 
for an absolute constant $\gamma >0$, which is more restrictive than the bound of Theorem 1.2. 

One reason why our inductive proof of Theorems 1.2 and 2.3 achieve less restrictive conditions than (4.1.4) is that we aim to prove that $p(z)$ has no zeros in some fixed, but otherwise arbitrarily small neighborhood of $[0, 1] \subset {\Bbb C}$, which is still sufficient to deduce Theorem 1.2 from, while the cluster expansion conditions ensure a stronger claim of $p(z)$ having no zeros in a disc.
In our case, for the disc (4.1.3) to contain a neighborhood of $[0, 1]$, we need to have $\rho=1+\delta$ for some fixed $\delta >0$.  If we modify our proof so that it works also for the disc (4.1.3), we would also arrive to (4.1.4),  basically because the condition (2.3.1) that $|f_i(x)| \leq 1 + {1 \over 6 \Delta}$ would have to be replaced by $|f_i(x)| \leq 2 + \delta$, which would require the values of ${\bold P}(A_i)$ to be much smaller. On the other hand, one advantage of having the polynomial $p(z)$ to be non-zero in the disc (4.1.3) containing $[0, 1]$ is that it allows one to improve the bound for $K$ in Theorem 1.2 to $K=O(\ln(n/\epsilon))$, where the implicit constant in the ``$O$'' notation is absolute.

In their approximation of the volume of the truncated hypergraph matching polytope, Guo and N proposed in \cite{GN25} a different polymer model, where the underlying graph is bipartite, with the coordinates on one side matched to the events on the other side. It would be interesting to find out if the model extends to the general case treated in this paper, and if it does, what bounds on the probabilities of $A_i$ it gives. 

\subhead (4.2) Optimality \endsubhead It is unclear if the upper bound for ${\bold P}(A_i)$ in Theorem 1.2 needs to be exponentially small in $\mu_i$. From the proof of Theorem 2.3, the exponential dependence arises because of the following phenomenon: even when an event $A \subset {\boldsymbol \Omega_1} \times {\boldsymbol \Omega_2}$ has small probability, for some $\omega_1 \in \Omega_1$, its section 
$$A_{\omega_1}=\bigl\{ \omega_2 \in {\boldsymbol \Omega_2}: \ (\omega_1, \omega_2) \in A \bigr\} \tag4.2.1$$ may still have a large probability in ${\boldsymbol \Omega_2}$.
Since by the Markov inequality, the probability to hit such an $\omega_1 \in {\boldsymbol \Omega_1}$ cannot be too large, our inductive proof carries through, but the very existence of such $\omega_1$ creates the exponential in $\mu_i$ dependence in our proof.

When we have a better control over the probability of $A_{\omega_1}$ defined by (4.2.1), the bounds for ${\bold P}(A_i)$ improve and may cease to be exponential in $\mu_i$.
This is the case, for example, when in Theorem 1.2 each set $A_i \subset {\boldsymbol \Omega}$ is a direct product 
$A_i = A_i^1 \times \cdots \times A_i^m$ with $A_i^j \subset \Omega_j$ or a union of a small number of such products. This situation occurs in counting solutions to CNF Boolean formulas, independent sets  in hypergraphs and in counting hypergraph colorings,  cf. \cite{BB25}, \cite{G+24}, \cite{HWY23}, \cite{JPV22}, \cite{L+25}, \cite{Moi19}.

On the other hand, the bounds obtained in \cite{GN25} for ``bad events" $A_i$ in the problem of approximating the volumes of some combinatorially defined polytopes indeed appear to be exponentially small in $\mu_i$. More precisely, the bounds of \cite{GN25}, when stated in terms of our Theorem 1.2, appear to be 
${\bold P}(A_i)=\Delta^{-O(r)} r^{-O(r)}$ for $r=\max_{i=1, \ldots, n} r_i$.


\Refs
\widestnumber\key{AAAAA}

\ref\key{Bar08}
\by A. Barvinok
\book Integer Points in Polyhedra
\bookinfo Zurich Lectures in Advanced Mathematics
\publ European Mathematical Society (EMS)
\publaddr  Z\"urich
\yr  2008
\endref

\ref\key{Bar16}
\by  A. Barvinok
\book Combinatorics and Complexity of Partition Functions
\bookinfo Algorithms and Combinatorics  {\bf 30}
\publ Springer
\publaddr Cham
\yr 2016
\endref

\ref\key{Bar17}
\by A. Barvinok
\paper Approximating permanents and hafnians
\jour Discrete Analysis 
\yr 2017
\pages Paper No. 2, 34 pp
\endref

\ref\key{BR19}
\by A. Barvinok and G. Regts
\paper Weighted counting of solutions to sparse systems of equations
\jour Combinatorics, Probability and Computing
\vol 28 
\yr 2019
\pages no. 5, 696--719
\endref

\ref\key{BB25}
\by F. Bencs and P. Buys
\paper Optimal zero-free regions for the independence polynomial of bounded degree hypergraphs
\jour Random Structures $\&$ Algorithms
\vol 66
\pages no. 4, Paper No. e70018, 32 pp
\yr 2025
\endref

\ref\key{BR24}
\by F. Bencs and G. Regts
\paper Approximating the volume of a truncated relaxation of the independence polytope
\paperinfo {\tt arXiv:2404.08577}  
\yr 2024
\endref

\ref\key{BR25}
\by F. Bencs and G. Regts
\paper Barvinok's interpolation method meets Weitz's correlation decay approach
\paperinfo preprint {\tt arXiv:2507.03135}
\yr 2025
\endref

\ref\key{EL75}
\by P. Erd\H{o}s and L.  Lov\'asz
\paper Problems and results on 3-chromatic hypergraphs and some related questions. Infinite and finite sets 
\inbook Colloq., Keszthely, 1973; dedicated to P. Erd\H{o}s on his 60th birthday), Vols. I, II, III
\bookinfo Colloquia Mathematica Societatis J\'anos Bolyai, Vol. 10
\pages 609--627
\publ North-Holland Publishing Co.
\publaddr Amsterdam-London
\yr 1975
\endref

\ref\key{GK94}
\by P. Gritzmann and V.  Klee
\paper On the complexity of some basic problems in computational convexity. II. Volume and mixed volumes
\inbook Polytopes: Abstract, Convex and Computational (Scarborough, ON, 1993)
\pages  373--466
\bookinfo NATO Advanced Science Institutes Series C: Mathematical and Physical Sciences
\publ  Kluwer Academic Publishers Group
\publaddr Dordrecht
\yr 1994
\endref

\ref\key{GN25}
\by H. Guo and V. N. 
\paper Deterministic approximation for the volume of the truncated fractional matching polytope
\inbook LIPIcs. Leibniz International Proceedings in Informatics
\vol 325
\publ Schloss Dagstuhl. Leibniz-Zentrum f\"ur Informatik
\publaddr Wadern
\yr 2025
\pages  Art. No. 57, 14 pp
\endref


\ref\key{G+24}
\by D. Galvin, G. McKinley, W. Perkins, M. Sarantis and P. Tetali
\paper On the zeroes of hypergraph independence polynomials
\jour Combinatorics,  Probability and  Computing 
\vol 33 
\yr 2024
\pages no. 1, 65--84
\endref

\ref\key{HWY23}
\by K. He, C. Wang and Y. Yin
\paper Deterministic counting Lov\'asz local lemma beyond linear programming
\inbook Proceedings of the 2023 Annual ACM-SIAM Symposium on Discrete Algorithms (SODA)
\pages  3388--3425
\publ Society for Industrial and Applied Mathematics (SIAM)
\publaddr Philadelphia, PA
\yr 2023
\endref


\ref\key{Jen24}
\by M. Jenssen
\paper The cluster expansion in combinatorics
\inbook Surveys in Combinatorics 2024
\bookinfo London Math. Soc. Lecture Note Ser.
\vol 493
\pages 55--88
\publ Cambridge University Press
\publaddr  Cambridge 
\yr 2024
\eds F. Fischer and R. Johnson
\endref

\ref\key{Jer24}
\by M. Jerrum
\paper Fundamentals of partial rejection sampling
\jour Probability Surveys 
\vol 21 
\yr 2024
\pages 171---199
\endref

\ref\key{JPV22}
\by V. Jain, H.T.  Pham, and T.D. Vuong
\paper Towards the sampling Lov\'asz Local Lemma
\inbook 2021 IEEE 62nd Annual Symposium on Foundations of Computer Science -- FOCS 2021
\pages 173--183
\publ IEEE Computer Society
\publaddr Los Alamitos, CA
\yr 2022
\paperinfo also preprint {\tt arXiv:2102.08342}
\endref

\ref\key{KLS96}
\by J. Kahn, N. Linial and A. Samorodnitsky
\paper Inclusion-exclusion: exact and approximate
\jour Combinatorica 
\vol 16 
\yr 1996
\pages  no. 4, 465--477
\endref

\ref\key{KP86}
\by R. Koteck\'y and D. Preiss
\paper Cluster expansion for abstract polymer models
\jour Communications in Mathematical Physics
\vol 103 
\yr 1986
\pages no. 3, 491--498
\endref

\ref\key{LN90}
\by N. Linial and N. Nisan
\paper Approximate inclusion-exclusion
\jour Combinatorica 
\vol 10 
\yr 1990
\pages no. 4, 349--365
\endref

\ref\key{LSS19}
\by J. Liu, A. Sinclair and P. Srivastava
\paper The Ising partition function: zeros and deterministic approximation
\jour Journal of Statistical Physics
\vol 174
\yr 2019
\pages  no.2, 287--315
\endref

\ref\key{L+25}
\by J. Liu, C. Wang, Y. Yin and Y. Yu
\paper Phase transitions via complex extensions of Markov chains
\inbook STOC '25: Proceedings of the 57th Annual ACM Symposium on Theory of Computing
\pages 903--914
\yr 2025
\publ ACM
\publaddr New York 
\paperinfo also preprint {\tt arXiv:2411.06857}
\endref

\ref\key{Moi19}
\by A. Moitra
\paper Approximate counting, the Lov\'asz local lemma, and inference in graphical models
\jour Journal of the  ACM 
\vol 66 
\yr 2019 
\pages no. 2, Art. 10, 25 pp
\endref

\ref\key{MT10}
\by R.A. Moser and G. Tardos
\paper A constructive proof of the general Lov\'asz local lemma
\jour Journal of the ACM 
\vol 57 
\yr 2010
\pages no. 2, Art. 11, 15 pp
\endref

\ref\key{PR17}
\by V. Patel and G. Regts
\paper Deterministic polynomial-time approximation algorithms for partition functions and graph polynomials
\jour SIAM Journal on Computing 
\vol 46 
\yr 2017
\pages no. 6, 1893--1919
\endref

\ref\key{SS05}
\by A.D. Scott and A.D. Sokal
\paper The repulsive lattice gas, the independent-set polynomial, and the Lov\'asz local lemma
\jour Journal of Statistical  Physics 
\vol 118 
\yr 2005
\pages  no. 5-6, 1151--1261
\endref

\ref\key{She85}
\by J.B. Shearer
\paper On a problem of Spencer
\jour Combinatorica 
\vol 5 
\yr 1985
\pages no. 3, 241--245
\endref

\ref\key{WY24}
\by C. Wang and Y. Yin
\paper A sampling Lov\'asz local lemma for large domain sizes
\inbook IEEE 65th Annual Symposium on Foundations of Computer Science -- FOCS 2024
\publ IEEE Computer Society
\publaddr Los Alamitos, CA
\yr 2024
\endref


\endRefs
\enddocument
\end